\documentclass[a4paper,12pt]{article}
\usepackage{cmap}                        % Поддержка поиска русских слов в PDF (pdflatex)
\usepackage[cp1251]{inputenc}            % Выбор языка и кодировки
\usepackage[english]{babel}
\usepackage[left=2cm,right=2cm,top=1.7cm,bottom=1.7cm]{geometry} % поля страницы
\usepackage[ruled,vlined]{algorithm2e}
\usepackage{amssymb}
\usepackage{amsmath, amsthm}
\theoremstyle{plain}
\newtheorem{thm}{Theorem}
\newtheorem{lem}{Lemma}
\newtheorem{prop}{Proposition}
\newtheorem{cor}{Corollary}

\theoremstyle{definition}

\newtheorem{pr}{Problem}

\begin{document}

\begin{center}\large
\textbf{A test for a local formation of finite groups to be a formation of  soluble groups with the Shemetkov property\footnote{This work is supported by the Ministry of Education of Belarus (``Convergence--2025'', 20211750)}}\normalsize

\smallskip
V.\,I. Murashka

 \{mvimath@yandex.ru\}

Department of  Mathematics and Technologies of Programming,

 Francisk Skorina Gomel State University, Gomel, Belarus\end{center}

\textbf{Abstract.}
L.A. Shemetkov posed a Problem 9.74    in Kourovka Notebook to find all local formations $\mathfrak{F}$ of finite groups such that every finite minimal non-$\mathfrak{F}$-group is either a Schmidt group  or a group of prime order. All known solutions to this problem are obtained under the assumption that every minimal non-$\mathfrak{F}$-group is soluble.
Using the above mentioned solutions we present a polynomial in $n$ time check for a local formation $\mathfrak{F}$ with bounded $\pi(\mathfrak{F})$ to be a formation of soluble groups with the Shemtkov property where $n=\max \pi(\mathfrak{F})$.



\section*{Introduction and the results}

All groups considered here are finite. Recall that a non-nilpotent group all whose proper subgroups are nilpotent is called a Schmidt group in the honor of O. Yu. Schmidt who described the structure of such groups \cite{Shmidt1924} in 1924. In 1951 N. Ito \cite[Proposition 2]{Ito1951} proved that a non-$p$-nilpotent group all whose proper subgroups are $p$-nilpotent is a Schmidt group.
Recall that for a class of groups $\mathfrak{X}$ a group $G\not\in\mathfrak{X}$ is called a minimal non-$\mathfrak{X}$-group if all its proper subgroups belong $\mathfrak{X}$. V.N. Semenchuk and A.F. Vasil’ev \cite{Semenchuck1984} described all hereditary local formations of soluble groups $\mathfrak{F}$ such that every soluble minimal non-$\mathfrak{F}$-group is either a Schmidt group  or a group of prime order. In 1984 L.A. Shemetkov asked  \cite[Problem 9.74]{Kour}    in Kourovka Notebook  to find all local formations $\mathfrak{F}$ of finite groups such that every finite minimal non-$\mathfrak{F}$-group is either a Schmidt group  or a group of prime order. The solutions to this problem are presented in \cite[Corollary 1]{Kamornikov1994} and \cite[Theorem 2]{BallesterBolinshes1995}.
Recall that a formation $\mathfrak{F}$ is said to have the Shemetkov property if every  minimal non-$\mathfrak{F}$-group is either a Schmidt group  or a group of prime order.

\begin{thm}[{\cite[Corollary 2.4.23]{Kamornikov2003}}]\label{thm0}
Let $\mathfrak{F}$ be a hereditary local formation. Then $\mathfrak{F}$ has the Shemetkov property if and only if it satisfies the following conditions:

$1)$ Every minimal non-$\mathfrak{F}$-group is soluble;

$2)$ $\mathfrak{F}$ is locally defined by $f$ where $f(p_i)=\mathfrak{G}_{\pi_i}$ for all $p_i\in\pi(\mathfrak{F})$ where $\pi_i$ is a subset of $\pi(\mathfrak{F})$ with $p_i\in\pi_i$.
\end{thm}

It is natural to ask if condition 1) of Theorem \ref{thm0} can be deduced from condition 2) of this theorem. The aim of this paper is to solve a particular case of this question: can one deduce from condition $2)$ of Theorem \ref{thm0} that $\mathfrak{F}$ is a formation of soluble groups with the Shemetkov property?

\begin{thm}\label{thm1}
 Let    $\pi=\{p_1, p_2,\dots, p_k\}$ be a set of primes not greater than $n$, $\pi_i$ be a subset of $\pi$ with $p_i\in\pi_i$. Assume that $\mathfrak{F}$ is a local formation with $\pi(\mathfrak{F})=\pi$ locally defined by $f$ where $f(p_i)=\mathfrak{G}_{\pi_i}$. In $O(n^2)$ operations one can check whether $\mathfrak{F}$ is a formation of soluble groups with the Shemetkov property.
 \end{thm}

In the section ``Proof of Corollary \ref{cor0}'' we will show how this algorithm works on a simple example.

\begin{cor}\label{cor0}
  Let $\pi=\{2, 3, 5, 7\}$ and $\mathfrak{F}$ be a local formation with $\pi(\mathfrak{F})=\pi$ locally defined by $f$ where $f(2)=f(3)=\mathfrak{G}_{\{2, 3, 5, 7\}}$, $f(5)=\mathfrak{G}_{\{3, 5, 7\}}$ and $f(7)=\mathfrak{G}_{\{5, 7\}}$. Then $\mathfrak{F}$ is a formation of soluble groups with the Shemetkov property.
\end{cor}

The proof of Theorem \ref{thm1} is based on the concept of $N$-critical graph. Recall that a Schmidt $(p, q)$-group is a Schmidt group with the normal Sylow $p$-subgroup. An \emph{$N$-critical graph} $\Gamma_{Nc}(G)$ of a group $G$ \cite[Definition 1.3]{VM} is a directed graph on the vertex set $\pi(G)$  and $(p, q)$ is an edge of   $\Gamma_{Nc}(G)$ iff $G$ has  a Schmidt $(p, q)$-subgroup. We can modify the proof of Theorem \ref{thm1} to show

\begin{cor}\label{cor1}
  Let $\Gamma$ be a directed graph such that $V(\Gamma)$ is a finite set of primes   and $n=\max V(\Gamma)$. One can check if every group $G$ with $\Gamma_{Nc}(G)=\Gamma$ is soluble in a polynomial in $n$ time.
\end{cor}

From the proves of Corollaries \ref{cor1} and \ref{cor2} follows

\begin{cor}\label{cor2}
 If a $\{2, 3, 5,7\}$-group $G$ does not contain  Schmidt $(5, 2)$-subgroups, $(7, 2)$-subgroups and $(7,3)$-subgroups, then $G$ is soluble.
\end{cor}

\section{Preliminaries}

All unexplained notations and terminologies are standard. The reader is referred to \cite{s9, Doerk1992} if necessary.
Here $Z_n$ is the cyclic group of order $n$;  $\pi(G)$ is the set of all prime divisors of  $|G|$;  $\pi(\mathfrak{X})=\underset{G\in\mathfrak{X}}\cup\pi(G)$;  $\Phi(G)$ is the Frattini subgroup of $G$; $\mathcal{M}(\mathfrak{X})$  is the class of all minimal non-$\mathfrak{X}$-groups.
Recall that $\mathfrak{G}_\pi$ is the class of all $\pi$-groups and $\mathbb{F}_p$ denotes a field with $p$ elements. Whenever
$V$ is a $G$-module over $\mathbb{F}_p$, $V\rtimes G$ denotes the semidirect product of $V$ with $G$
corresponding to the action of $G$ on $V$ as $G$-module.

\subsection{Formations}

Recall that a \emph{formation} is a class of groups $\mathfrak{F}$ which is   closed  under taking epimorphic images (i.e. from $G\in\mathfrak{F}$ and $N\trianglelefteq G$ it follows that $G/N\in\mathfrak{F}$)  and subdirect products (i.e. from $G/N_1\in\mathfrak{F}$ and $G/N_2\in\mathfrak{F}$ it follows that $G/(N_1\cap N_2)\in\mathfrak{F}$).

A formation $\mathfrak{F}$ is   called

$(a)$ \emph{hereditary} if  from $G\in\mathfrak{F}$ and $H\leq G$ it follows that $H\in\mathfrak{F}$.

$(b)$ \emph{saturated} if  from $G/N\in\mathfrak{F}$ where $N\trianglelefteq G$  and $N\leq\Phi(G)$ it follows that $G\in\mathfrak{F}$.

A function of the form $f: \mathbb{P}\rightarrow\{formations\}$ is called a \emph{formation function}. Recall \cite[IV, Definitions 3.1]{Doerk1992} that a formation $\mathfrak{F}$ is called \emph{local} if $$\mathfrak{F}=(G\mid G/C_G(H/K)\in f(p)\textrm{ for every } p\in\pi(H/K) \textrm{ and every chief factor }H/K\textrm{ of }G)$$ for some formation function $f$. In this case $f$ is called a \emph{local definition} of $\mathfrak{F}$.

\subsection{$N$-critical graph}

Here a (directed) graph $\Gamma$  is a pair of sets $V(\Gamma)$ and $E(\Gamma)$ where $V(\Gamma)$ is a set of vertices of $\Gamma$ and
$E(\Gamma)$ is a set of edges of $ \Gamma$, i.e. the set of ordered pairs of elements from $V(\Gamma)$.  Two graphs $\Gamma_1$ and $\Gamma_2$ are called equal (denoted by $\Gamma_1 = \Gamma_2$) if $V (\Gamma_1) = V (\Gamma_2)$ and $E(\Gamma_1) = E(\Gamma_2)$. Graph
$\Gamma_1$ is called a subgraph of $\Gamma_2$ (denoted by $\Gamma_1\subseteq\Gamma_2$) if $V (\Gamma_1) \subseteq V (\Gamma_2)$ and $E(\Gamma_1) \subseteq E(\Gamma_2)$. Graph $\Gamma$ is called
a union of graphs $\Gamma_1$ and $\Gamma_2$ (denoted by $\Gamma = \Gamma_1\cup\Gamma_2)$ if $V(\Gamma) = V (\Gamma_1)\cup V (\Gamma_2)$ and $E(\Gamma) = E(\Gamma_1)\cup E(\Gamma_2)$.

The $N$-critical graph of a class of groups $\mathfrak{X}$ \cite[Definition 3.1]{VM} is defined by
 $$\Gamma_{Nc}(\mathfrak{X})=\bigcup_{G\in\mathfrak{X}}\Gamma_{Nc}(G).$$

\begin{lem}[{\cite[Theorem 2.7(2)]{VM}}]\label{lem1}
  Let $H$   be a subgroup of $G$  and $G_i$ be groups where $1\leq i\leq n$. Then $\Gamma_{Nc}(H)\subseteq\Gamma_{Nc}(G)$ and $$\Gamma_{Nc}\left(\times_{i=1}^n G_i\right)=\bigcup_{i=1}^n\Gamma_{Nc}(G_i).$$\end{lem}

\begin{prop}[{\cite[Proposition 6.1]{VM}}]\label{prop1} The following statements hold:

  $(a)$ If $p$ is a prime, then
$V(\Gamma_{Nc}(PSL(2, 2^p)))=\pi(2(2^{2p}-1))$  and
$E(\Gamma_{Nc}(PSL(2, 2^p)))=\{(2, q) \mid q\in\pi(2^p-1)\}\cup \{(q, 2) \mid q\in\pi(2^{2p}-1)\}.$

$(b)$ If $p$ is an odd prime, then  $V(\Gamma_{Nc}(PSL(2, 3^p)))=\pi(3(3^{2p}-1))$ and
 $E(\Gamma_{Nc}(PSL(2, 3^p)))=\{(3, q) \mid q\in\pi(3^p-1)\setminus\{2\}\}\cup \{(2, 3)\}\cup \{(q, 2) \mid q\in\pi(3^{2p}-1)\setminus\{2\}\}$.

$(c)$ If $p>5$ is a prime with  $p^2+1\equiv 0 \mod 5$, then  $V(\Gamma_{Nc}(PSL(2, p)))=\pi(p(p^{2}-1))$ and
 $E(\Gamma_{Nc}(PSL(2, p)))=\{(p, q) \mid q\in\pi\left(\frac{p-1}{2}\right)\}\cup \{(2, 3)\}\cup \{(q, 2) \mid q\in\pi(p^2-1)\setminus\{2\}\}$.

$(d)$ If   $p$ is an odd prime, then   $V(\Gamma_{Nc}(Sz(2^p)))=\pi(2(2^{2p}+1)(2^p-1))$ and
$E(\Gamma_{Nc}(Sz(2^p)))=\{(2, q) \mid q\in\pi(2^p-1)\}\cup \{(q, 2) \mid q\in\pi((2^p-1)(2^{2p}+1))\}$.

$(e)$  $V(\Gamma_{Nc}(PSL(3,3)))=\{2, 3, 13\}$  and $E(\Gamma_{Nc}(PSL(3,3)))=\{(2, 3), (3, 2), (13, 3)\}$.

\end{prop}

\subsection{Algorithms}

We assume that basic arithmetic operations and comparison are done in the same amount of time (equal 1). Recall that   $O(f(n))$ means a function $g(n)$ such that there exist $C>0$ and a natural number $n_0$ such that for all $n\geq n_0$ holds $|g(n)|<C|f(n)|$. It is well known that the decomposition of a natural number $n$ into the product of primes can be done in $O(\sqrt{n})$ operations.
\section{Proof of Theorem \ref{thm1}}

$(a)$ \emph{$\mathfrak{F}$ is hereditary formation}.

Note that $f(p)$ is a hereditary formation for all $p\in\pi=\pi(\mathfrak{F})$. Hence $\mathfrak{F}$ is a hereditary formation by
\cite[IV, Proposition 3.14]{Doerk1992}.

$(b)$ \emph{$\mathfrak{F}$ is a formation of soluble groups with the Shemetkov property if and only if $\Gamma_{Nc}(G)\not\subseteq \Gamma_{Nc}(\mathfrak{F})$ for every minimal simple non-abelian group $G$}.

Suppose that $\Gamma_{Nc}(G)\not\subseteq \Gamma_{Nc}(\mathfrak{F})$ for every minimal simple non-abelian group $G$.
 Assume  that $\mathfrak{F}$ contains a non-soluble group $G_1$. Since $\mathfrak{F}$ is a hereditary, it contains a minimal non-soluble group $G_2$. Note that $G_3\simeq G_2/\Phi(G_2)\in\mathfrak{F}$ is  a minimal simple non-abelian group. Hence
 $\Gamma_{Nc}(G_3)\subseteq\Gamma_{Nc}(\mathfrak{F})$, a contradiction.

 Thus every group in $\mathfrak{F}$ is soluble. Assume now that $\mathcal{M}(\mathfrak{F})$ contains a non-soluble group $G_1$. Note that $G_1$ is a minimal non-soluble group. By Gasch\"{u}tz-Lubeseder-Schmid Theorem $\mathfrak{F}$ is a saturated formation. Hence a minimal simple non-abelian group $G_2\simeq G_1/\Phi(G_1)\in\mathcal{M}(\mathfrak{F})$. Since $G_2$ is not a Schmidt group or a group of prime order, $\Gamma_{Nc}(G_2)$ is the join of $\Gamma_{Nc}(M)$ where $M$ runs through all maximal subgroups of $G_2$. From $M\in\mathfrak{F}$ for every maximal subgroup $M$ of $G_2$ it follows that $\Gamma_{Nc}(G_2)\subseteq\Gamma_{Nc}(\mathfrak{F})$, a contradiction. Thus every group in $\mathcal{M}(\mathfrak{F})$ is soluble. Therefore $\mathfrak{F}$ has the Shemetkov property by Theorem \ref{thm0} or \cite[Theorem 6.4.12]{s9}.

Suppose now that $\mathfrak{F}$ is a formation of soluble groups with the Shemetkov property. Assume  that there is a  minimal simple non-abelian group $G$ with $\Gamma_{Nc}(G)\subseteq\Gamma_{Nc}(\mathfrak{F})$. Since $G$ is non-soluble, it contains a minimal non-$\mathfrak{F}$-group $H$ as a subgroup. Since $\pi(G)=V(\Gamma_{Nc}(G))\subseteq V(\Gamma_{Nc}(\mathfrak{F}))$ and $\mathfrak{F}$ has the Shemetkov property, we see that $H$ is a Schmidt $(p, q)$-group for some $(p, q)\not\in E(\Gamma_{Nc}(\mathfrak{F}))$. From $\Gamma_{Nc}(H)\subseteq\Gamma_{Nc}(G)$ by Lemma \ref{lem1} it follows that
$\Gamma_{Nc}(G)\not\subseteq\Gamma_{Nc}(\mathfrak{F})$, a contradiction. Thus $\Gamma_{Nc}(G)\not\subseteq \Gamma_{Nc}(\mathfrak{F})$ for every minimal simple non-abelian group $G$.

$(c)$ $E(\Gamma_{Nc}(\mathfrak{F}))=\{(p_i, p_j)\mid p_i\in\pi, p_j\in\pi_i\}$.

Let $\Gamma$ be a graph with $V(\Gamma)=\pi$ and $E(\Gamma)=\{(p_i, p_j)\mid p_i\in\pi, p_j\in\pi_i\}$.
Assume that $(p, q)\in E(\Gamma_{Nc}(\mathfrak{F}))$. Note that $p=p_i$ and $q=p_j$ for some $p_i, p_j\in\pi$ with $p_i\neq p_j$. Then $\mathfrak{F}$ contains a Schmidt $(p_i, p_j)$-group $G_1$.
Hence $G_2\simeq G_1/\Phi(G_1)\in\mathfrak{F}$ is also a Schmidt $(p_i, p_j)$-group. In this case $G$ has the unique minimal normal subgroup $N$, $N$ is a $p_i$-group and $G/C_G(N)\simeq Z_{p_j}$.  Therefore $Z_{p_j}\in f(p_i)=\mathfrak{G}_{\pi_i}$. Hence $p_j\in\pi_i$. Thus $E(\Gamma_{Nc}(\mathfrak{F}))\subseteq E(\Gamma)$.

Let $p_i\in \pi$ and $p_j\in\pi_i$.
According to \cite[B, Theorem 10.3]{Doerk1992} a group $Z_{p_j}$ has a faithful irreducible module $V$ over a field $\mathbb{F}_{p_i}$.
Let $G=V\rtimes Z_{p_j}$. Note that $V$ is the unique minimal normal subgroup of $G$ and every maximal subgroup of $G$ is either a $p_i$-group or a $p_j$-group. Hence $G$ is a Schmidt $(p_i, p_j)$-group. From $G/C_G(V)\simeq Z_{p_j}\in f(p_i)$ and $f(p_j)\neq \emptyset$ it follows that $G\in\mathfrak{F}$. Hence $(p_i, p_j)\in E(\Gamma_{Nc}(\mathfrak{F}))$.
Therefore $E(\Gamma)\subseteq E(\Gamma_{Nc}(\mathfrak{F}))$. Thus $E(\Gamma)=E(\Gamma_{Nc}(\mathfrak{F}))$.

$(d)$  \emph{Let $p$ be a prime and $p\leq n$. If $G$ is isomorphic to any group from $PSL(2, 2^p)$, $PSL(2, 3^p)$ and
$ Sz(2^p)$ for odd $p$, $PSL(2, p)$ for $p>5$   with $5\in\pi(p^2+1)$, then we can check if $\Gamma_{Nc}(G)\subseteq\Gamma$ in $O(n)$ operations}.

First we need to check that $\pi(G)\subseteq \pi$. Note that $$|G|\in\{(2^{p}+1)(2^{2p}-2^p), (3^{p}+1)(3^{2p}-3^p)/2, 2^{2p}(2^{2p}+1)(2^p-1), (p+1)(p^2-p)/2\}.$$ Since $p\leq n$, the rude estimates shows that $|G|\leq 2^{6n}$. Hence $|G|$ has no more than $6n$ not necessary different primes divisors. Note that if $\Gamma_{Nc}(G)\subseteq\Gamma$, then all these divisors belong to $\pi$. Since $|\pi|<n$, we see that in no more than $7n$ divisions we can check the if $\pi(G)\subseteq\pi=V(\Gamma)$.
Suppose now that $\pi(G)\subseteq\pi$.

Assume that $G\simeq PSL(2, 2^p)$. Since $\pi(G)\subseteq \pi$, in $O(n)$ operations we can compute $\pi(2^p-1)$ and $\pi(2^p+1)$. Note that $|\pi(2^p-1)\cup\pi(2^p+1)|<n$ in our case. According to Proposition \ref{prop1}(a)  taking vertex 2, we need to check that from it starts an arrow to every vertex from $\pi(2^p-1)$; taking every vertex from $\pi(2^p-1)\cup\pi(2^p+1)$, we need to check that from it starts an arrow to   vertex 2. It is clear that all these can be done in $O(n)$ operations.

Assume that $G\simeq Sz(2^p)$. Since $\pi(G)\subseteq \pi$, in $O(n)$ operations we can compute $\pi(2^p-1)$ and $\pi(2^{2p}+1)$. Note that $|\pi(2^p-1)\cup\pi(2^{2p}+1)|<n$ in our case. According to Proposition \ref{prop1}(d) taking vertex 2, we need to check that from it starts an arrow to every vertex from $\pi(2^p-1)$; taking every vertex from $\pi(2^p-1)\cup\pi(2^{2p}+1)$, we need to check that from it starts an arrow to   vertex 2. It is clear that all these can be done in $O(n)$ operations.

Assume that $G\simeq PSL(2, 3^p)$. Since $\pi(G)\subseteq \pi$, in $O(n)$ operations we can compute $\pi(3^p-1)$ and $\pi(3^p+1)$. Note that $|\pi(3^p-1)\cup\pi(3^p+1)|<n$ in our case. According to Proposition \ref{prop1}(b) taking vertex 3, we need to check that from it starts an arrow to every vertex from $\pi(3^p-1)\setminus\{2\}$; taking every  vertex from $\pi(3^p-1)\cup\pi(3^p+1)\setminus\{2\}$, we need to check that from it starts an arrow to   vertex 2; also we need to check, if $(2, 3)\in E(\Gamma)$.   It is clear that all these can be done in $O(n)$ operations.

Assume that $G\simeq PSL(2, p)$. Since $\pi(G)\subseteq \pi$, in $O(n)$ operations we can compute $\pi(p-1)$ and $\pi(p+1)$. Note that $|\pi(p-1)\cup\pi(p+1)|<n$ in our case. According to Proposition \ref{prop1}(c) taking vertex $p$, we need to check that from it starts an arrow to every vertex from $\pi(\frac{p-1}{2})$; taking every  vertex from $\pi(p-1)\cup\pi(p+1)\setminus\{2\}$, we need to check that from it starts an arrow to   vertex 2; also we need to check, if $(2, 3)\in E(\Gamma)$.   It is clear that all these can be done in $O(n)$ operations.

$(e)$ \emph{In $O(n^{3/2})$ operations we can show that there are up to isomorphism no more than $4n$ minimal simple non-abelian groups $G$ for which $\Gamma_{Nc}(G)\subseteq \Gamma$ is possible and list this groups}.

Recall \cite[II, Bemerkung 7.5]{Huppert1967} that   minimal simple non-abelian groups up to isomorphism are $PSL(2, 2^p)$ for a prime $p$, $PSL(2, 3^p)$ and $Sz(2^p)$ for an odd prime $p$, $PSL(2, p)$ where $p>5$ is a prime with $5\in\pi(p^2+1)$ and $PSL(3,3)$.

  %$(b.1)$ \emph{If $G\simeq PSL(2, 2^p)$, or $G\simeq Sz(2^p)$ or , then $p\in\pi(q-1)$ for some $q\in\pi$}.

Assume that $G\simeq PSL(2, 2^p)$ or $G\simeq Sz(2^p)$. Then $\pi(2^p-1)\subseteq \pi$. Hence if $q\in\pi(2^p-1)$, then $2^p\equiv 1\mod q$. Note that $2^{q-1}\equiv 1\mod q$ by Fermat's Little Theorem. Hence $2^{(q-1,p)}\equiv 1\mod q$. Since $p$ is a prime, $p\in\pi(q-1)$.

Assume that $G\simeq PSL(2, 3^p)$ where $p$ is an odd prime. Then $\pi(3^p-1)\setminus\{2\}\subseteq \pi$. So there is  $q\in\pi(3^p-1)\setminus \{2, 3\}$. By analogy $p\in\pi(q-1)$.

Let $\rho=\cup_{q\in\pi}(\pi(q-1))$. Note $\pi(q-1)$ can be computed in $O(n^{1/2})$ divisions for every $q\in\pi$ and $\max\rho<n$. From $|\pi|<n$ it follows that  $\rho$ can be computed in $O(n^{3/2})$ operations. Hence there are up to isomorphism no more than $3|\rho|<3n$  minimal simple groups $G$  from $$\{PSL(2, 2^p), PSL(2, 3^p),  Sz(2^p)\mid p\in\mathbb{P}\}$$ for which $\Gamma_{Nc}(G)\subseteq \Gamma$ is possible. Note that if $\Gamma(PSL(2, p))\subseteq\Gamma$, then $p\in\pi$.
Since $|\pi|\leq n$,  in $O(n^{3/2})$ operations we can show that there are up to isomorphism no more than $4n$ minimal simple non-abelian groups $G$ for which $\Gamma_{Nc}(G)\subseteq \Gamma$ is possible and list this groups.

$(f)$ \emph{The final step}.

According to $(b)$ we need only to check that $\Gamma_{Nc}(G)\not\subseteq \Gamma_{Nc}(\mathfrak{F})$ for every minimal simple non-abelian group $G$. Using $(e)$ we can list (up to isomorphism) no more than $4n$ such groups $G$ for which
$\Gamma_{Nc}(G)\subseteq \Gamma_{Nc}(\mathfrak{F})$ is possible in $O(n^{3/2})$ operations. According to $(d)$ for every listed group $G$ in $O(n)$ operations we can check if $\Gamma_{Nc}(G)\subseteq \Gamma_{Nc}(\mathfrak{F})$ (note that we can check if $\Gamma_{Nc}(PSL(3,3))\subseteq \Gamma_{Nc}(\mathfrak{F})$ in $O(1)$ operations).  Thus in
$O(n^{2})$ operations we can check if  $\Gamma_{Nc}(G)\not\subseteq \Gamma_{Nc}(\mathfrak{F})$ for every minimal simple non-abelian group $G$.

\begin{algorithm}[H]
\caption{IsSolubleSFormation($\mathfrak{F}$)}
\SetAlgoLined
\KwResult{True, if $\mathfrak{F}$ is a formation of soluble groups with the Shemtkov property and false otherwise.}
\KwData{$\pi=\{p_1, p_2,\dots, p_k\}$ is a set of primes not greater than $n$, $\pi_i$ is a subset of $\pi$ with $p_i\in\pi_i$ for $i\in\{1,\dots,k\}$}

 $\Gamma\gets$ graph  with $V(\Gamma)=\pi$ and $(p_i, p_j)\in E(\Gamma)$ iff $p_j\in\pi_i$\;
%n^2

\If{$\{(2, 3), (3, 2), (13, 3)\}\subseteq E(\Gamma)$ \textbf{or} $\{(2, 3), (3, 2), (5, 2)\}\subseteq E(\Gamma)$}{\Return false\;}

\For{$p\in\pi$ with $5\in\pi(p^2+1)$ and $p>5$}
   {\If{$\pi((p^3-p)/2)\subseteq\pi$}{
   \If{$ \{(p, q) \mid q\in\pi\left(\frac{p-1}{2}\right)\}\cup \{(2, 3)\}\cup \{(q, 2) \mid q\in\pi(p^2-1)\setminus\{2\}\}  \subseteq E(\Gamma)$}
   {\Return false\;}
   }}

$\rho\gets\cup_{p\in\pi}(\pi(p-1)\setminus\{2\})$\;

\For{$p\in\rho$}
  {\If{$\pi(2(2^{2p}-1))\subseteq\pi$}{
     \If{$ \{(2, q) \mid q\in\pi(2^p-1)\}\cup \{(q, 2) \mid q\in\pi(2^{2p}-1)\}  \subseteq E(\Gamma)$}
   {\Return false\;}}
   \If{$\pi(2(2^p-1)(2^{2p}+1))\subseteq\pi$}{
 \If{$ \{(2, q) \mid q\in\pi(2^p-1)\}\cup \{(q, 2) \mid q\in\pi((2^p-1)(2^{2p}+1))\}  \subseteq E(\Gamma)$}
   {\Return false\;}}
 \If{$\pi(3(3^{2p}-1)/2)\subseteq\pi$}{
 \If{$ \{(3, q) \mid q\in\pi(3^p-1)\setminus\{2\}\}\cup \{(2, 3)\}\cup \{(q, 2) \mid q\in\pi(3^{2p}-1)\setminus\{2\}\}  \subseteq E(\Gamma)$}
   {\Return false\;}
  }}

 \Return{true\;}
 \end{algorithm}

\section{Proof of Corollary \ref{cor0}}

Let do the steps according to ``IsSolubleSFormation'':\\
1.   $V(\Gamma_{Nc}(\mathfrak{F}))=\{2, 3, 5, 7\}$ and
 $$E(\Gamma_{Nc}(\mathfrak{F}))=\{(2, 3), (2, 5), (2, 7), (3, 2), (3, 5), (3, 7), (5, 3), (5, 7), (7,5)\}.$$
2. $\{(2, 3), (3, 2), (13, 3)\}\not\subseteq E(\Gamma)$ and $\{(2, 3), (3, 2), (5, 2)\}\not\subseteq E(\Gamma)$.\\
3. The  only prime $p>5$ in $\pi$ with $5\in\pi(p^2+1)$ is 7. Note that $\{(7, 3), (2, 3), (3, 2)\}\not\subseteq E(\Gamma)$.\\
4. $\rho=(\pi(2-1)\cup\pi(3-1)\cup\pi(5-1)\cup\pi(7-1))\setminus\{2\}=\{3\}$.\\
5. $\pi(2(2^6-1))\subseteq \pi$ but $\{(2, 7), (7, 2), (3, 2)\}\not\subseteq E(\Gamma)$.\\
6. $\pi(2(2^3-1)(2^6+1))\not\subseteq \pi$.\\
7. $\pi(3(3^6-1)/2)\not\subseteq \pi$\\
8. Thus $\mathfrak{F}$ is a formation of soluble groups with the Shemetkov property.

\section{Proof of Corollary \ref{cor1}}

Let $\pi=\{p_1,\dots, p_k\}=V(\Gamma)$, $\pi_i=\{p_i\}\cup \{p_j\mid (p_i, p_j)\in E(\Gamma)\}$ and $\mathfrak{F}$ be a local formation with $\pi(\mathfrak{F})=\pi$ locally defined by $f$ where $f(p_i)=\mathfrak{G}_{\pi_i}$ for all $p_i\in\pi$. Since $n=\max \pi$, in a polynomial in $n$ time we can check does $\mathfrak{F}$ is a formation of soluble groups with the Shemetkov property by Theorem \ref{thm1}.

Let prove that every group $G$ with $\Gamma_{Nc}(G)=\Gamma$ is soluble if and only if $\mathfrak{F}$ is a formation of soluble groups with the Shemetkov property.  Note that $\Gamma_{Nc}(\mathfrak{F})=\Gamma$ by step $(c)$ of the proof of Theorem \ref{thm1}.

Suppose that every group $G$ with $\Gamma_{Nc}(G)=\Gamma$ is soluble. Assume that there is a minimal simple non-abelian group $G_1$ with $\Gamma_{Nc}(G_1)\subseteq\Gamma$. Let $H(p_i, p_j)$ be a Schmidt $(p_i, p_j)$-group and $$G_2=G_1\times (\times_{(p_i,p_j)\in E(\Gamma)}H(p_i, p_j))\times (\times_{p_i\in V(\Gamma)}Z_{p_i}).$$ Since $E(\Gamma)\cup V(\Gamma)$ is finite, by Lemma \ref{lem1}
$$\Gamma_{Nc}(G_2)=\Gamma_{Nc}(G_1)\cup \left(\bigcup_{(p_i,p_j)\in E(\Gamma)}\Gamma_{Nc}(H(p_i, p_j))\right)\cup \left(\bigcup_{p_i\in V(\Gamma)}\Gamma_{Nc}(Z_{p_i})\right)=\Gamma.$$
Therefore there exists a non-soluble group $G_2$ with $\Gamma_{Nc}(G_2)=\Gamma$, a contradiction. Thus $\Gamma_{Nc}(G_1)\not\subseteq\Gamma$ for every minimal simple non-abelian group $G_1$. Thus $\mathfrak{F}$ is a formation of soluble groups with the Shemetkov property by  step $(b)$ of the proof of Theorem \ref{thm1}.

Suppose that $\mathfrak{F}$ is a formation of soluble groups with the Shemetkov property. Assume that there exists a non-soluble group $G$ with $\Gamma_{Nc}(G)=\Gamma$. Then $G$ contains a minimal non-$\mathfrak{F}$-group $H$ as a subgroup. Since $\pi(G)=\pi=\pi(\mathfrak{F})$, $H$ is a Schmidt $(p, q)$-group for some $(p, q)\not\in E(\Gamma)$. Hence $\Gamma_{Nc}(G)\neq\Gamma$, a contradiction. Thus every group $G$ with $\Gamma_{Nc}(G)=\Gamma$ is soluble.

\section{Final remarks}

Let $\mathfrak{F}, \pi, \pi_i$ be the same as in Theorem \ref{thm1}. Then if ``IsSolubleSFormation($\mathfrak{F}$)'' returns true, then we automatically prove a lot of results about $\mathfrak{F}$. Lets list some of them.

From \cite{Kazarin1991} and \cite{Vasilev1993} or \cite{BallesterBolinches2005} it follows that $\mathfrak{F}$ has \emph{the Kegel property}, i.e.
\begin{center}
  If $A,B,C\in\mathfrak{F}$ and $G=AB=BC=CA$, then $G\in \mathfrak{F}$.
\end{center}
From \cite[Lemma 2.2]{Amberg1998} it follows that $\mathfrak{F}$ has \emph{the property $\mathcal{P}_2$}, i.e.
\begin{center}
  If $G=A_1\dots A_k$ where $A_iA_j\in\mathfrak{F}$ for $i\neq j$, then $G\in\mathfrak{F}$.
\end{center}
From    \cite{Vasilev1987} (see also \cite[Theorem 1]{Vasilev2022}) it follows that $\mathfrak{F}$ has \emph{the Belonogov property} in the class of all soluble groups
\begin{center}
   If   $A, B, C$ are non-conjugate maximal $\mathfrak{F}$- subgroups of a soluble group $G$, then $G\in \mathfrak{F}$.
\end{center}

It is well known that the class of all nilpotent groups can be characterized as a class of groups whose every Sylow (or cyclic primary) subgroup is subnormal. Recall \cite[Definition 6.1.2]{s9} that a subgroup $H$ of a group $G$ is said to be $\mathfrak{F}$-subnormal in $G$ if either $H = G $ or there exists a chain of subgroups
$H =H_0 <\dots< H_n = G$ such that $H_{i-1}$ is a  maximal subgroup of $H_i$ and $H_i/Core_{H_i}(H_{i-1})\in\mathfrak{F}$ for $i=1,\dots, n$. From \cite[Corollaries 3.9 and 3.10]{Murashka2020} it follows that for $\mathfrak{F}$ hold the analogues of this characterization.
\begin{center}
  If every Sylow (or every cyclic primary) subgroup of $G$ is $\mathfrak{F}$-subnormal in $G$, then  $G\in\mathfrak{F}$.
\end{center}

From \cite[Theorem 4.4]{VM} or \cite[Theorem 1]{BallesterBolinches1996} it follows that for $\mathfrak{F}$ holds an analogue of Frobenious $p$-nilpotency criterion.

\begin{center}
A $\pi$-group $G$ belongs $\mathfrak{F}$ iff $N_G(P)/C_G(P)$ is a $\pi_i$-group for every $p_i$-subgroup\,$P$\,of\,$G$\,and\,$p_i\in\pi$.
\end{center}

We want to note that the following problems seems interesting and remains open.

\begin{pr}
  Let    $\pi=\{p_1, p_2,\dots, p_k\}$ be a set of primes not greater than $n$, $\pi_i$ be a subset of $\pi$ with $p_i\in\pi_i$. Assume that $\mathfrak{F}$ is a local formation with $\pi(\mathfrak{F})=\pi$ locally defined by $f$ where $f(p_i)=\mathfrak{G}_{\pi_i}$. Can one check that $\mathfrak{F}$ has the Shemetkov property in polynomial in $n$ time?
\end{pr}

{\small\bibliographystyle{siam}
\bibliography{Graph7}}

\end{document}